\documentclass[fleqn,10pt,a4paper]{article}
\usepackage{style}
\usepackage[hypertex]{hyperref}

\DeclareMathOperator{\B}{Bar}
\DeclareMathOperator{\Ext}{Ext}
\DeclareMathOperator{\Tor}{Tor}
\DeclareMathOperator{\Hom}{Hom}
\DeclareMathOperator{\co}{co}
\def\HH{H\!H}
\let\to\rightarrow
\let\ot\otimes
\let\lact\rightharpoonup
\let\ract\leftharpoonup
\let\eps\varepsilon
\let\*\bullet

\def\cf{\textit{cf.}}

\def\loccit{\textit{loc.~cit.}}

\makeatletter
\def\operator@font{\sf}
\makeatother

\title{On the cohomology of a Galois entwining}

\author{Mariano Suarez Alvarez\thanks{This work was supported by a grant from UBACYT TW62 and the
PICT~03-08280 project.} 
}

\date{October 4th, 2004}

\begin{document}

\maketitle

\section{Introduction}

\paragraph An extension of algebras $A/B$ is said to be Hopf-Galois over a Hopf algebra $H$ if $A$
is an $H$-comodule algebra, $B=A^{\co H}$ is the subalgebra of coinvariants, and a certain Galois
condition is satisfied, \cf~\cite{M}. Such an extension can be viewed as a non-commutative principal
bundle with structure group $H$. It turns out, though, that to develop a satisfactory theory of
principal bundles in the non-commutative setting the notion of Hopf-Galois extensions is too
restrictive; it does not apply, for example, to all the quantum spheres of Podle\'s~\cite{P}.

One considers then a more general situation in which the r\^ole of the structure group is
played by a coalgebra. T.~Brzezi\'nski and P. Hajac~\cite{Br2} have
proposed a corresponding notion of \emph{coalgebra Galois extensions}, tightly related to that of
entwining structures introduced in~\cite{Br3}. In this context, one can fit the Podle\'s spheres in 
coalgebra Galois extensions $SU_q(2)/S_{q,s}^2$ for appropriate coactions of coalgebras on
$SU_q(2)$

\paragraph In a later paper~\cite{Br1}, Brzezi\'nski introduced two cohomology theories for an
entwining structure, and, in particular, for coalgebra $C$-Galois extensions $A/B$: the entwined
cohomology $H_\psi^\*(A,-)$ of $A$ with values in $A$-bimodules, and a $C$-equivariant version. He 
computed the entwined cohomology of $A$ when the algebra $B$ of coinvariants in $A$ is the ground field
$k$, and noted it is essentially trivial.

\paragraph The purpose of the present note is to record the extension of Brze\-zi\'nski's computation
of cohomology to the general case of a flat coalgebra Galois extension $A/B$. We show below that in that
situation $H_\psi(A,-)$ coincides with the Hochschild cohomology $\HH^\*(B,-)$ of the subalgebra $B$
of coinvariants.

We plan to study the equivariant comohology of the entwined structure corresponding to a coalgebra
Galois extension in a future paper.

\section{Coalgebra Galois extensions and the theorem}

\paragraph Fix a field $k$. All spaces and (co)algebras considered below are $k$-vector spaces and
$k$-(co)algebras, and all unadorned tensor products are taken over $k$. Most of our statements
can be extended to the slightly more general situation in which $k$ is simply a ring, provided one
adds appropriate projectivity or flatness hypotheses.

\paragraph An \textit{entwining structure} is a triple $(A,C,\psi)$ consisting of an algebra $A$, a
coalgebra $C$ and a map $\psi:C\ot A\to A\ot C$\dash---which we write \`a la
Sweedler, with implicit sums over greek indices, as in $\psi(c\ot a)=a_\alpha\ot
c^\alpha$\dash---satifying the following compatibility conditions:
  \begin{align*}
  &(aa')_\alpha\ot c^\alpha=a_\alpha a'_\beta\ot c^{\alpha\beta}, 
  &&1_\alpha\ot c^\alpha = 1\ot c, \\
  &a_\alpha\ot{c^\alpha}_1\ot{c^\alpha}_2=a_{\beta\alpha}\ot{c_1}^\alpha\ot{c_2}^\beta, 
  &&a_\alpha\eps(c^\alpha) = a\eps(c).
  \end{align*}

\paragraph\label{p:gal} Entwining structures arise naturally in the following situation. Let $C$ be a coalgebra.
Let $A$ be an algebra which is a right $C$-comodule, and let 
  \[
  B=\{b\in A:(ba)_0\ot(ba)_1=ba_0\ot a_1 \text{for all $a\in A$}\}.
  \]
This is a subalgebra. There is a linear map $\beta:A\ot_BA\to A\ot C$ such that $\beta(a\ot
a')=aa'_0\ot a'_1$, which is evidently left $A$-linear and right $C$-colinear; when $\beta$ is
bijective, we say that $A/B$ is a \emph{$C$-Galois extension}. If this is the case, we let
$\gamma:C\to A\ot_BA$ be the unique map such that $\beta\circ\gamma=\eta\ot1$, which we shall write
$\gamma(c)=l(c)\ot r(c)$ with an implicit summation over an implicit index. Then there is a
canonical entwining structure $(A,C,\psi)$ associated to the extension $A/B$ in which the map $\psi$
is given by $\psi(c\ot a) = \beta(\gamma(c)a)=l(c)(r(c)a)_0\ot(r(c)a)_1$.

\paragraph Let $(A,C,\psi)$ be an entwining structure. We shall always consider the space $A\ot C$
to be endowed with the structure of an $A$-bimodule with left and right actions given by
  \begin{align*}
  &\lambda\lact a\ot c=\lambda a\ot c, &&a\ot c\ract\rho=a\lact\psi(c\ot\rho)
  \end{align*}
for each $a,\lambda,\rho\in A$, $c\in C$.

\paragraph We note that when $A/B$ is a $C$-Galois extension, the Galois map $\beta:A\ot_BA\to A\ot C$ 
is a map of $A$-bimodules. Indeed, we have
  \begin{align*}
  \beta(a\ot a')\ract b 
    &= aa'_0\ot a'_1\ract b
    = aa'_0l(a'_1)(r(a'_1)b)_0\ot(r(a'_1)b)_1 \\
    &= a(a'b)_0\ot(a'b)_1 
    = \beta(a\ot a'\ract b),
  \end{align*}
where the third equality follows from the fact, stated as property \textit{(iii)} in the proof of theorem 2.7
in \cite{Br2}, that $a_0l(a_1)\ot r(a_1)=1\ot a$ for all $a\in A$. 

\paragraph\label{p:res} In \cite{Br1}, Brzezi\'nski considers the complex $\B^\psi_\*(A)=(A\ot
C)\ot_A\B_\*(A)$; here $\B_\*(A)$ is the usual Hochschild resolution of $A$ as an
$A$-bimodule. Since of course $A$ is flat as a left $A$-module, this complex is acyclic over $A\ot
C$, and since its components are clearly free as $A$-bimodules, we have in fact a projective
resolution of $A\ot C$ as an $A$-bimodule.

\paragraph For each $A$-bimodule $M$, \cite{Br1} defines the cohomology of the entwining structure $(A,C,\psi)$
with values in $M$ to be the graded space $H_\psi^\*(A,M)$ obtained by taking the homology of the
cochain complex $\Hom_{A^e}(\B^\psi_\*(A),M)$. In view of the observation made in \pref{p:res}, we
have at once that $H_\psi^\*(A,M)=\Ext_{A^e}^\*(A\ot C,M)$.

Observe that with this identification in mind, proposition 2.3 in \cite{Br1}, stating that $A\ot C$
is a projective $A$-bimodule iff $H_\psi^1(A,-)$ vanishes identically, becomes immediate.

\paragraph Proposition 2.6 in \cite{Br1} and the comments after its proof hint that when $A/B$ is a
$C$-Galois extension, the cohomology of the corresponding entwining structure $(A,C,\psi)$ is
related to the Hochschild cohomology of $B$. In that paper the case where $B=k$ is considered; we
have, more generally,

\begin{Theorem}
Let $C$ be a coalgebra. Let $A/B$ be a $C$-Galois extension, and let $(A,C,\psi)$ be the
corresponding entwining structure. Then we have $H_\psi^0(A,-)\cong H^0(B,-)$ as functors of $A$-bimodules.
In fact, if $A$ is flat as a (left or right) $B$-module, $H_\psi^\*(A,-)\cong H^\*(B,-)$ as
$\partial$-functors on the category of $A$-bimodules. 
\end{Theorem}

\begin{Proof}
Because $\beta:A\ot_BA\to A\ot C$ is an isomorphism of $A$-bimodules,
  \begin{equation}\label{eq:iso}
  H_\psi^\*(A,-)\cong\Ext_{A^e}^\*(A\ot C,-)\cong\Ext_{A^e}^\*(A\ot_BA,-)
  \end{equation}
naturally on $A$-bimodules. On the other hand, the change-of-rings spectral sequence \textsc{XVI}.\S5.$(2)_3$
constructed in \cite{CE}, when specialised to the morphism $B^e\to A^e$, has
$E_2^{p,q}\cong\Ext_{A^e}^q(\Tor^{B^e}_p(A^e,B),-)$ and converges to $\Ext_{B^e}^\*(B,-)$; note that we know
from corollary \textsc{IX}.\S4.4, \loccit, that $\Tor^{B^e}_\*(A^e,B)\cong\Tor^B_\*(A,A)$.
Now, since this spectral sequence lives on the first quadrant, \eqref{eq:iso} and convergence
immediately imply that $H_\psi^0(A,-)\cong H^0(B,-)$. When $A$ is flat as a $B$-module, the spectral
sequence degenerates at once, and this, together with \eqref{eq:iso}, gives an isomorphism
$\Ext_{A^e}^\*(A\ot_BA,-)\cong\Ext_{B^e}^\*(B,-)$.~\qed 

\end{Proof}


\vfill
\begin{center}
\begin{minipage}{0.9\textwidth}
\small
\noindent
Departamento de Matem\'atica.\\
Facultad de Ciencias Exactas y Naturales. 
Universidad de Buenos Aires.\\
Ciudad Universitaria. Pabell\'on I.
Buenos Aires (1428) Argentina.
\par\noindent\textit{E-mail address:} \texttt{mariano@dm.uba.ar}
\par\noindent\textit{URL:} \texttt{http://mate.dm.uba.ar/\~{}aldoc9}
\end{minipage}
\end{center}

\end{document}